\documentclass[12pt]{amsart}

\usepackage{amssymb}
\usepackage{amsmath}
\usepackage{latexsym}
\usepackage{graphics}
\usepackage{color}

\newcommand{\A}{\ensuremath{\mathbb{A}}}
\newcommand{\Q}{\ensuremath{\mathbb{Q}}}
\newcommand{\Z}{\ensuremath{\mathbb{Z}}}

\newcommand{\F}{\ensuremath{\mathbb{F}}}

\newcommand{\Zp}{\ensuremath{{\mathbb Z}_p}}

\DeclareMathOperator{\Frac}{Frac}

\DeclareMathOperator{\mres}{mres}
\DeclareMathOperator{\LRings}{\cL(Rings)}
\DeclareMathOperator{\LMR}{\cL(\cM\cR)}
\DeclareMathOperator{\LRingsMR}{\cL(Rings, \cM\cR)}

\DeclareMathOperator{\Spec}{Spec}

\DeclareMathOperator{\Dim}{Dim}

\newcommand\cL{\mathcal{L}}
\newcommand\cM{\mathcal{M}}

\newcommand\cO{\mathcal{O}}

\newcommand\cR{\mathcal{R}}

\newcommand\frm{\mathfrak{m}}

\newtheorem{theorem}[subsection]{Theorem}

\newtheorem{lemma}[subsection]{Lemma}

\newtheorem{remark}[subsection]{Remark}

\newtheorem{namedTheorem}[subsection]{\theoremname}
\newcommand{\theoremname}{testing}
\newenvironment{named}[1]{\renewcommand\theoremname{#1}
\begin{namedTheorem}}
{\end{namedTheorem}}
\begin{document}

\title{Geometric proofs of theorems of Ax-Kochen and Er{\v{s}}ov}

\author[Denef]{Jan Denef}

\dedicatory{Dedicated to the memory of Professor Jun-ichi Igusa}

\address[]{KULeuven, Department of Mathematics, Celestijnenlaan 200 B, B-3001 Leuven (Heverlee), Belgium.}
\email{jan.denef@wis.kuleuven.be}

\thanks{We thank Dan Abramovich, Steven Dale Cutkosky, and Kalle Karu for stimulating conversations and information.}

%\date{\today}
\date{October 30, 2015}

\begin{abstract}
We give an algebraic geometric proof of the Theorem of Ax and Kochen on p-adic diophantine equations in many variables. Unlike Ax-Kochen's proof, ours does not use any notions from mathematical logic and is based on weak  toroidalization of morphisms. We also show how this geometric approach yields new proofs of the Ax-Kochen-Er{\v{s}}ov transfer principle for local fields, and of quantifier elimination theorems of Basarab and Pas.

\end{abstract}

\maketitle
%\tableofcontents

\section{Introduction}

The purpose of the present paper is to give purely algebraic geometric proofs of the  following theorem of Ax and Kochen \cite{AxKoch-1} and of some other related results, such as quantifier elimination, which go back to Ax, Kochen, and Er{\v{s}}ov. 

\begin{theorem} \textbf{\emph{Ax-Kochen's Theorem on p-adic forms.}} For any positive integer $d$ there exists a positive integer $N$ such that for each prime $p > N$, any homogeneous polynomial over the ring $\Z_p$ of p-adic integers, with degree $d$ and more than $d^2$ variables, has a nontrivial zero over $\Zp$.
\label{TheoremAxKochenForms}
\end{theorem}

The proof of Ax-Kochen, and also ours, starts from the well known elementary fact that the above theorem is true for $\Z_p$ replaced by the ring $\F_p[[t]]$ of formal power series over the field $\F_p$ with $p$ elements (in this case one can take $N=1$).
Thus, in order to prove Theorem \ref{TheoremAxKochenForms}, it suffices to prove the next theorem (applying it to the universal family of hypersurfaces of degree $d$ in projective space of dimension $d^2$).

\begin{theorem} \textbf{\emph{Transfer of surjectivity.}}
\label{TheoremTransferSurjectivity}
Let $f: X \rightarrow Y$ be a morphism of integral separable schemes of finite type over $\Z$.
For all large enough primes $p$ we have:
$$f: X(\Z_p) \rightarrow Y(\Z_p)$$
is surjective if and only if
$$f: X(\F_p[[t]]) \rightarrow Y(\F_p[[t]])$$
is surjective.
\end{theorem}

In the present paper we prove Theorem \ref{TheoremTransferSurjectivity} by using the Theorem of Abramovich and Karu \cite{Abramovich-Karu} on Weak Toroidalization of Morphisms (extended to non-closed fields in \cite{Abramovich-Denef-Karu}). Instead of weak toroidalization  one could use Cutkosky's Theorem on Local Monomialization \cite{Cutkosky-LocMon}, see Remark \ref{RemarkAboutTameness2}. In \cite{Denef-ColliotConjecture} we give a second algebraic geometric proof of Theorem \ref{TheoremAxKochenForms}, not using transfer to $\F_p[[t]]$, by proving a conjecture of Colliot-Th\'{e}l\`{e}ne \cite{Colliot-1}. Our proof \cite{Denef-ColliotConjecture} of this conjecture is also based on weak toroidalization of morphisms.

Ax and Kochen obtained Theorem \ref{TheoremAxKochenForms} as a direct consequence of the following more powerful Transfer Principle due to Ax and Kochen \cite{AxKoch-1} and (independently) Er{\v{s}}ov \cite{Ersov-1,Ershov-2}.

\begin{theorem} \textbf{\emph{Ax-Kochen-Er{\v{s}}ov Transfer Principle.}}
\label{TheoremAKE-TransferPrinciple}
Let $\varphi$ be an assertion  in the language $\LRings$ of rings (see section \ref{SubsectionLanguages} ).
For all large enough primes $p$ we have the following. The assertion $\varphi$ is true in $\Z_p$ if and only if it is true in $\F_p[[t]]$.
\end{theorem}
Ax-Kochen and Er{\v{s}}ov proved this Transfer Principle using methods from mathematical logic (model theory). An elementary but very ingenious proof has been
given by Cohen \cite{Cohen-QuantifierElimination}(see also Weispfenning \cite{Weispfenning-1}, and Pas \cite{Pas-1}), but his method is still very much in the spirit of mathematical logic. In the present paper we also give a new proof of this Transfer Principle, again based on weak toroidalization.

Moreover, in section \ref{SectionEliminationQuantifiers} we give a geometric proof of Basarab's Quantifier Elimination Theorem \cite{Basarab-1} for henselian valuation rings of residue field characteristic large enough or zero. Basarab's result is a refinement of a quantifier elimination theorem of Ax-Kochen \cite{AxKochen-3}, and is related to work of Delon \cite{Delon} and Weispfenning \cite{Weispfenning-2}. His proof uses the methods of Ax-Kochen and Er{\v{s}}ov, and is based on model theory. Basarab's theorem directly implies the quantifier elimination theorem of Pas \cite{Pas-1}, which has several applications in arithmetic algebraic geometry. It enables to study certain integrals over local fields  \cite{Denef-Degree},\cite{Pas-1},\cite{Cluckers-Gordon-1}, in particular generalizations of Igusa's local zeta functions \cite{Igusa-2,Igusa-Book-2}, and has several applications to motivic integration \cite{Denef-Loeser-GermsArcs},\cite{Cluckers-Loeser-1,Cluckers-Loeser-2,Cluckers-Loeser-3}. This relates to work of Lichtin \cite{Lichtin-Monomialization1,Lichtin-Monomialization-2}, who was the first to apply monomialization of morphisms (i.e. local toroidalization) to study multivariate Igusa fiber integrals. 

The present paper is about henselian valuation rings of residue field characteristic large enough or zero. 
Using multiplicative residues of higher order (i.e. with respect to certain proper ideals), instead of the multiplicative residues introduced in section \ref{SectionMultRes}, the method of the present paper can be easily adapted to give geometric proofs of quantifier elimination results of Basarab \cite{Basarab-1} and Pas \cite{Pas-2} that are valid for henselian valuation rings of characteristic zero,  without any restriction on the residue field.  This approach can be much simplified in case of $\Z_p$, for a fixed prime $p$, using compactness, to get an easy proof, based on weak toroidalization, of Macintyre's Quantifier Elimination Theorem \cite{Macintyre-pAdicElimination}; this is done in \cite{Denef-MonomializAndpAdicElimination}.

Our paper is organized as follows. In section \ref{SectionMultRes} we discuss multiplicative residues of elements in local integral domains. These also play a key role in Basarab's paper \cite{Basarab-1} on quantifier elimination.
In section \ref{SectionLogHensel} we formulate and prove a logarithmic version of Hensel's Lemma. We did not fit it in the framework of logarithmic geometry, but this has been done more recently by Cao \cite{Cao-ColliotConjecture}. 
We recall the Weak Toroidalization Theorem in section \ref{SectionToroidalizationMorphisms}. The heart of the present paper is section \ref{SectionTamenessTheorem} where we state and prove what we call the Tameness Theorem \ref{TheoremTameness}. The Weak Toroidalization Theorem reduces its proof to the case of a log-smooth morphism where the Tameness Theorem is a direct consequence of the logarithmic Hensel's Lemma. 
We prove Theorem \ref{TheoremTransferSurjectivity} on transfer of surjectivity in section \ref{SectionTransferSurjectivity} as an easy consequence of the Tameness Theorem and Lemma \ref{LemmaTransferResidues} on transfer of residues. This lemma is stated and proved in section \ref{SectionTransferResidues} as an easy application of embedded resolution of singularities.
Finally, in section \ref{SectionEliminationQuantifiers} we formulate and prove Basarab's Quantifier Elimination Theorem \ref{TheoremBasarab} using the Tameness Theorem, and we prove the Ax-Kochen-Er{\v{s}}ov Transfer Principle \ref{TheoremAKE-TransferPrinciple} as a direct consequence of Basarab's Theorem.

Our motivation to develop an algebraic geometric approach to
quantifier elimination for henselian valuation rings, comes from the above mentioned applications to the study of variants of the local zeta functions that Igusa introduced in \cite{Igusa-2,Igusa-Book-2}. We are happy to dedicate this work to the memory of late Professor Jun-ichi Igusa.

\subsection{Terminology and notation}
In the present paper, $R$ will always denote a noetherian integral domain. A \emph{variety over} $R$ is an integral separated scheme of finite type over $R$. A rational function $x$ on a variety $X$ over $R$ is called \emph{regular at} a point $P \in X$ if it belongs to the local ring $\cO_{X,P}$ of $X$ at $P$, and it is called \emph{regular} if it is regular at each point of $X$.

\emph{Uniformizing parameters over} $R$ on a variety $X$ over $R$, are regular rational functions on $X$ that induce an \'{e}tale morphism to an affine space over $R$.

A \emph{reduced strict normal crossings divisor over} $R$ on a smooth variety $X$ over $R$ is a closed subset $D$ of $X$ such for any $P \in X$ there exist uniformizing parameters $x_1, \dots, x_n$ over $R$ on an open neighborhood of $P$, such that for any irreducible component $C$ of $D$, containing $P$, there is an $i \in \{1, \dots, n\}$ which generates the ideal of $C$ in $\cO_{X,P}$.

% TODO: Toevoegen dat log-smooth impliceert toroidal ???

\section {Multiplicative residues} \label{SectionMultRes}

Let $R$ be a noetherian integral domain, and $X$ a variety over $R$.
Let $A$ be any local $R$-algebra which is an integral domain. We denote by $\frm_A$ its maximal ideal, by $\Frac(A)$ its field of fractions, and by $\eta_A$ the generic point of $\Spec(A)$. For any $A$-rational point $a  \in X(A)$ on $X$ we denote by $a \!\mod \frm_A$ the $A / \frm_A$-rational point on $X$ induced by $a$. For any $x \in \cO_{X,a(\eta_A})$ the pullback $a^*(x)$ of $x$ to $\Frac(A)$ is denoted by $x(a) \in \Frac(A)$. Moreover, for $a, a' \in X(A)$ we write $a \equiv a' \!\mod \frm_A$ to say that  $a \!\mod \frm_A \;=\; a' \!\mod \frm_A$.

\begin{named}{Definition} \label{DefinitionMultRes1} \rm
Let $z,z' \in \Frac(A)$. The elements $z,z'$ have \emph{same multiplicative residue} if
$$ z' \in z(1 + \frm_A).$$
\newline
Let $a,a' \in X(A)$ and let $x_1, \dots , x_r$ be rational functions on $X$. The points $a,a'$ have the \emph{same residues
with respect to} $x_1, \dots , x_r$ if $a \equiv a'\!\mod\frm_A$ and, for $i = 1, \dots, r$, the following two conditions hold.
\begin{enumerate}
\item The rational function $x_i$ is regular at $a(\eta_A)$ if and only if it is regular at  $a'(\eta_A)$.
\item $x_i(a), x_i(a') \in \Frac(A)$ have same  multiplicative residue if $x_i$ is regular at both $a(\eta_A)$ and $a'(\eta_A)$.
\end{enumerate}
\end{named}

Instead of working with rational functions $x_1, \dots , x_r$ we can also work with locally principal closed subschemes of $X$, i.e. subschemes whose sheaf of ideals is locally generated by a single element.
\begin{named}{Definition} \rm
Let $a,a' \in X(A)$ and $D_1, \dots , D_r$ locally principal closed subschemes of $X$. The points $a,a'$ have the \emph{same residues with respect to} $D_1, \dots , D_r$ if $a \equiv a'\!\mod\frm_A$ and, for $i = 1, \dots, r$, the following condition holds. 
Let $g_i \in \cO_{X,a(\frm_A)}$ be a generator for the ideal of $D_i$ at $a(\frm_A)$, then $g_i(a), g_i(a') \in A$ have same multiplicative residue.
\end{named}

\begin{lemma}
\label{LemmaResidues1}
Let $X$ be an affine variety over $R$, and let $x_1, \dots, x_r$ be rational functions on $X$. Then there exist regular rational functions $x'_1, \dots, x'_{r'}$ on $X$ such that for any local R-algebra $A$, which is an integral domain, and any $a, a' \in X(A)$ we have the following. The points $a$ and $a'$ have the same residues with respect to $x_1, \dots, x_r$ if they have the same residues with respect to  $x'_1, \dots, x'_{r'}$.
\end{lemma}
\emph{Proof}. This is clear, by taking for $x'_1, \dots, x'_{r'}$ any finite list of regular rational functions on $X$ which satisfies the following condition. For each $i \in \{1, \dots, r\}$ and each $P \in X$ with $x_i$ regular at $P$, there are elements $x'_j$ and $x'_k$ in this list with $x_i = x'_j/x'_k$, and $x'_k(P)\neq 0$.  Obviously, such a finite list exists if $X$ is affine. $\square$

\begin{lemma}
\label{LemmaResidues2}
Assume that $R$ is a noetherian normal integral domain. Let $X$ be a smooth variety over $R$, and $x_1, \dots, x_r$  regular rational functions on $X$. Let $D \subset X$ be a reduced strict normal crossings divisor over $R$ containing the zero locus of each $x_i$, for $i= 1, \dots, r$. Then, for any local R-algebra $A$, which is an integral domain, and any $a, a' \in X(A)$ we have the following. The points $a$ and $a'$ have the same residues with respect to $x_1, \dots, x_r$ if they have the same residues with respect to the irreducible components of $D$.
\end{lemma}

\emph{Proof}. Assume that $a, a' \in X(A)$ have the same residues with respect to the irreducible components of $D$. Set $P := a(\frm_A)$. Then also $P = a'(\frm_A)$. Let $D_1, \dots, D_m$ be the irreducible components of $D$ that contain $P$. Since $D$ has normal crossings over $R$, there exist uniformizing parameters $z_1, \dots, z_n$ over $R$ on an open neighborhood of $P$, with $n \geq m$, such that $z_i$ generates the ideal of $D_i$ in $\cO_{X,P}$ for $i=1, \dots, m$. 
We can write each $x_i$, for $i=1, \dots, r$, as a monomial in $z_1, \dots, z_m$ times a unit $u$ in $\cO_{X,P}$, because $X$ is normal (being smooth over a normal ring). Obviously, this implies the lemma, because $u(a)$ and $u(a')$ have same multiplicative residue if $a \equiv a'\!\mod\frm_A$. $\square$

\subsection{The structure $\cM\cR(A)$ of multiplicative residues} \label{SubSecStructureMR}

Assume now that $R = \Z$, thus $A$ is any local integral domain. Denote by $\cM\cR(A) := A/(1+\frm_A)$ the set of equivalence classes of the equivalence relation $\{(z,z')\in A\times A \; | \, z' \in z(1+\frm_A)\}$ on $A$. The equivalence class of an element $a \in A$ is called the \emph{multiplicative residue of} $a$ and is denoted by $\mres(a)$. Note that $\cM\cR(A)$ is a commutative monoid with identity $1:=\mres(1)$ and with multiplication induced by the one on $A$. Moreover $\cM\cR(A)$ is equipped with the natural multiplicative map
$s_A: \cM\cR(A) \rightarrow A / \frm_A: \mres(a) \mapsto a + \frm_A$ onto the residue field $A / \frm_A$. This map  induces an isomorphism from the group $\cM\cR(A)^\times$ of elements of $\cM\cR(A)$ which have an inverse, onto the multiplicative group of the residue field. Indeed the elements of $\cM\cR(A)^\times$ are precisely the multiplicative residues of the units of $A$. Note that  $0:=\mres(0)$ is the unique element of the monoid  $\cM\cR(A)$ that multiplied with any element of this monoid equals itself.
We denote by  $+_{\mathrm{mod}}$ the  binary composition law on $\cM\cR(A)$ defined as follows. For $\alpha, \beta \in \cM\cR(A)$, the composition $\alpha +_{\mathrm{mod}} \beta$ is the unique $\gamma$ in $\cM\cR(A)^\times \cup \{0\}$ with $s_A(\gamma) = s_A(\alpha) + s_A(\beta)$. 

Let $B$ be a any local integral domain. We call a bijection  $\tau: \cM\cR(A) \rightarrow \cM\cR(B)$ an \emph{isomorphism} if it is compatible with multiplication and if there exists a (necessarily unique) field isomorphism $\overline{\tau}: A / \frm_A \rightarrow B / \frm_B$ such that $\overline{\tau} \circ s_A = s_B \circ \tau$. These conditions are equivalent to the requirement that $\tau$ is compatible with multiplication and $+_{\mathrm{mod}}$.
If $A$ and $B$ are valuation rings, then obviously such an isomorphism $\tau$ also induces an isomorphism of ordered groups between the value groups of $\Frac(A)$ and $\Frac(B)$.

\section{Logarithmic Hensel's Lemma} \label{SectionLogHensel}

Let $R$ be a noetherian normal integral domain. Let $f : X \rightarrow Y$ be a morphism of smooth varieties over $R$, and $D \subset X$, $E \subset Y$ reduced strict normal crossings divisors over $R$ with $f^{-1}(E) \subset D$. Denote the irreducible components of $D$ and $E$ by $D_1, \dots, D_r$ and $E_1, \dots, E_s$.
Note that $X$ and $Y$ are normal, because they are smooth over the normal ring $R$.

\subsection{Log-smoothness} 
Let $P \in X$.
Choose uniformizing parameters $x_1, \dots , x_n$ over $R$ on an open neighborhood of $P$ in $X$ so that locally at $P$ the locus of $\prod_i x_i$ is $D$.
Choose uniformizing parameters $y_1, \dots , y_m$ over $R$ on an open neighborhood of $f(P)$ in $Y$ so that locally at $f(P)$ the locus of $\prod_j y_j$ is $E$.
Recall that \emph{uniformizing parameters over} $R$ on a variety $U$ over $R$, are regular rational functions on $U$ that induce an \'{e}tale morphism to an affine space over $R$.
Since $R$ is normal, each such uniformizing parameter generates a prime ideal or is a unit in the local ring of any point of $U$.

\textbf{Definition}
The morphism $f$ is called \emph{log-smooth at }$P$ with respect to $D, E$, if the logarithmic jacobian
$$\left( \frac{\partial \log f^*(y_j)}{\partial \log x_i}(P)\right )_{i=1, \dots, n,\, j=1, \dots, m} $$
(considered as a matrix over the residue field of $P$) has rank equal to the relative dimension of $Y/R$ at $f(P)$.
\newline Here as usual $\frac{\partial \log f^*(y_j)}{\partial \log x_i}$ denotes the rational function $\frac{x_i}{f^*(y_j)} \frac{\partial f^*(y_j)}{\partial x_i}$ on $X$.
\newline Note that $\frac{\partial \log f^*(y_j)}{\partial \log x_i}$ belongs to $\cO_{X,P}$, because $ f^*(y_j)$ can be written as a unit in  $\cO_{X,P}$ times a monomial in $x_1, \dots, x_n$, since $f^{-1}(E) \subset D$ and $X$ is normal. Note also that the above definition of log-smoothness does not depend on the choice of the uniformizing parameters $x_i, y_j$.

\begin{named}{Logarithmic Hensel's Lemma}
Let $A$ be any henselian local $R$-algebra which is an integral domain, and $\frm_A$ its maximal ideal.
Let $a_0 \in X(A) \setminus D(A)$. Assume that f is log-smooth at $a_0(\frm_A) \in X$ with respect to $D, E$.
Then, any $b \in Y(A)$ having the same residues with respect to $E_1, \dots, E_s$ as $f(a_0)$, is the image under $f$ of an $a \in X(A)$ with the same residues as $a_0$ with respect to $D_1, \dots, D_r$.
\end{named}

\emph{Proof.} By restricting to suitable open neighborhoods of $a_0(\frm_A)$ and $f(a_0(\frm_A))$ we can make the following two assumptions.
\begin{description}
\item[a]
There exist uniformizing parameters $x_1, \dots , x_n$ over $R$ on $X$ such that the locus of $\prod_i x_i$ is $D$, and 
such that the locus of each $x_i$ is irreducible or empty.
\item[b]
There exist uniformizing parameters $y_1, \dots , y_m$ over $R$ on $Y$ such that the locus of $\prod_j y_j$ is $E$, and 
such that the locus of each $y_j$ is irreducible or empty.
\end{description}

\subsubsection{Changing coordinates}
Let $p: X\otimes_R A \rightarrow \A^n_A$ be the \'{e}tale morphism induced by $x_1, \dots , x_n$, and let $q: Y\otimes_R A \rightarrow \A^m_A$ be the \'{e}tale morphism induced by $y_1, \dots , y_n$. Consider the morphisms
\begin{align*}
&\alpha: \;  \A^n_A \rightarrow \A^n_A : (X_1, \dots , X_n) \mapsto (x_1(a_0)X_1, \dots ,x_n(a_0) X_n),
\\
&\beta: \; \A^m_A \rightarrow \A^m_A : (Y_1, \dots , Y_m) \mapsto (y_1(f(a_0))Y_1, \dots ,y_m(f(a_0)) Y_m).
\end{align*}
Set $X_A = X \otimes_R A$, $Y_A = Y \otimes_R A$, and denote by $f_A$ the morphism $f_A: \; X_A \rightarrow Y_A$ induced by $f$. Denote by $\alpha_X$ and $\beta_Y$ the morphisms obtained by base change of $\alpha$ and $\beta$ as defined by the following two cartesian  squares:
$$
\begin{array}{lcr}

\begin{array}{lcl} 
	X' & \stackrel{p'}{\longrightarrow} & \A^n_A \\
    \downarrow \alpha_X & & \downarrow \alpha \\ 
    X_A & \stackrel{p}{\longrightarrow} &  \A^n_A
\end{array}
& \qquad &
\begin{array}{lcl} 
	Y' & \stackrel{q'}{\longrightarrow} & \A^m_A \\
    \downarrow \beta_Y & & \downarrow \beta \\ 
    Y_A & \stackrel{q}{\longrightarrow} &  \A^m_A
\end{array}

\end{array}
$$

We think of $\alpha_X$ and $\beta_Y$ as coordinate changes induced by $x_i \mapsto x_i(a) x_i$ and $y_j \mapsto y_j(f(a)) y_j$, although these are not open immersions.

Denote the pull-back to $X'$ (through $p'$) of the standard affine coordinates on $\A^n_A$ by $x'_1, \dots , x'_n \in \Gamma(X', \cO_{X'})$. These are uniformizing parameters over $A$ on $X'$. 
Similarly, denote the pull-back to $Y'$ (through $q'$) of the standard affine coordinates on $\A^m_A$ by $y'_1, \dots , y'_m \in \Gamma(Y', \cO_{Y'})$. These are uniformizing parameters over $A$ on $Y'$.
By construction we have
\begin{equation} \label{equation-loghensel-1}
\alpha_X^*(x_i) = x_i(a_0)x'_i\, , \quad  
\beta_Y^*(y_j) = y_j(f(a_0))y'_j . 
\end{equation}

Note that $a_0 \in X(A)$ uniquely lifts to a point $a_0' \in X'(A)$ with $\alpha_X(a_0') = a_0$. Moreover $b \in Y(A)$ uniquely lifts to a point $b' \in Y'(A)$ with $\beta_Y(b') = b$, because $b$ has the same residues as $f(a_0)$ with respect to $y_1, \dots , y_m$. By \eqref{equation-loghensel-1} we have
\begin{equation} \label{equation-loghensel-2}
x_i'(a_0') = 1 \, , \quad  
y_j'(b') \equiv 1 \! \mod \frm_A. 
\end{equation}

\subsubsection{Lifting the morphism $f$}
We claim that there exists a morphism $f'$ of schemes over $A$ such that the following diagram commutes:
$$
\begin{array}{lcl} 
	X' & \stackrel{f'}{\longrightarrow} & Y' \\
    \downarrow \alpha_X & & \downarrow \beta_Y \\ 
    X_A & \stackrel{f_A}{\longrightarrow} &  Y_A
\end{array}
$$
To prove this we only have to show for $j=1, \dots, m$ that $\alpha_X^*(f_A^*(y_j))$ is divisible by $y_j(f(a_0))$ in $\Gamma(X', \cO_{X'})$.

By assumptions \textbf{a} and \textbf{b}, and because $X$ is normal and $f^{-1}(E) \subset D$, there exist non-negative integers $e_{i,j}$, for $i = 1, \dots, n$ and $j = 1, \dots, m$, and units $u_j \in \Gamma(X, \cO_{X})$ such that 
\begin{equation} \label{equation-loghensel-3}
f^*(y_j)=u_j\prod_{i=1}^n x_i^{e_{i,j}},
\end{equation}
for $j=1 \dots, m$. Hence $\alpha_X^*(f_A^*(y_j))$ is divisible by $\prod_{i=1}^n \alpha_X^*(x_i)^{e_{i,j}}$ in $\Gamma(X', \cO_{X'})$. Thus by (\ref{equation-loghensel-1}), $\alpha_X^*(f_A^*(y_j))$ is divisible by $\prod_{i=1}^n x_i(a_0)^{e_{i,j}}$ in $\Gamma(X', \cO_{X'})$.

On the other hand, evaluating(\ref{equation-loghensel-3}) on the rational point $a_0$ yields
$$
 y_j(f(a_0))=u_j(a_0)\prod_{i=1}^n x_i(a_0)^{e_{i,j}} \in A .
$$
Since $u_j(a_0)$ is a unit in $A$, we get that $y_j(f(a_0))$ divides $\prod_{i=1}^n x_i(a_0)^{e_{i,j}}$. Thus $y_j(f(a_0))$ divides $\alpha_X^*(f_A^*(y_j))$ in $\Gamma(X', \cO_{X'})$. This proves the claim.

\subsubsection{Applying the classical Hensel's Lemma}

Note that $f'(a_0')  \equiv b' \! \mod \frm_A$, because $f(a_0) \equiv b \! \mod \frm_A$ and 
$$
	y'_j(f'(a_0')) = 1 \equiv y'_j(b') \! \mod \frm_A,
$$
by the second equality in (\ref{equation-loghensel-1}) and the congruence in (\ref{equation-loghensel-2}).

Note also that the morphism $f'$ is smooth at 
$a'_0( \frm_A) \in X'$.
Indeed this follows from the log-smoothness of $f$ at $a_0(\frm_A)$, because, by the equations in (\ref{equation-loghensel-1}), the jacobian matrix of $f'$ at $a'_0(\frm_A)$, with respect to the uniformizing parameters $x'_1, \dots, x'_n$ and $y'_1, \dots, y'_m$, equals the logarithmic jacobian of $f$ at $a_0(\frm_A)$ 
with respect to $x_1, \dots, x_n, y_1, \dots, y_m$.

Hence, by the classical Hensel's Lemma for smooth morphisms, $b' \in Y'(A)$ can be lifted to a point $a' \in X'(A)$ with $f'(a') = b'$ and $a' \equiv a'_0 \! \mod \frm_A$.

\subsubsection{Conclusion of the proof of Logarithmic Hensel's Lemma.} Put $a := \alpha_X(a') \in X(A)$. Then $f(a) = b$, and 
$$
 a = \alpha_X(a') \equiv \alpha_X(a'_0) = a_0 \! \mod \frm_A.
$$
Moreover by (\ref{equation-loghensel-1}) we have
$$
 x_i(a) = \alpha_X^*(x_i)(a') = x_i(a_0) \, x'_i(a'). 
$$
But $x'_i(a') \equiv 1 \! \mod \frm_A$ by (\ref{equation-loghensel-2}), since $x'_i(a') \equiv x'_i(a'_0) \! \mod \frm_A$. Hence $x_i(a)$ and $x_i(a_0)$ have the same multiplicative residue. This finishes the proof  of the Logarithmic Hensel's Lemma. $\square$

\section{Toroidalization of morphisms} \label{SectionToroidalizationMorphisms}

\begin{named}{Definition} \rm
Let $K$ be a field of characteristic zero. Let $f: X\to Y$ be a dominant morphism of nonsingular varieties over $K$, and $D \subset X$, $E \subset Y$ reduced strict normal crossings divisors over $K$.
\newline
We call $f$ \emph{toroidal with respect to} $D$ and $E$ if $f^{-1}(E) \subset D$, and if, after base change to an algebraic closure $\overline{K}$ of $K$, for each closed point $a$ of $X \otimes_K \overline{K}$ there exist uniformizing parameters $x_1, \dots , x_n$ for $\widehat{\cO}_{X \otimes_K \overline{K}, \, a}$ and uniformizing parameters $y_1, \dots , y_m$ for $\widehat{\cO}_{Y \otimes_K \overline{K}, \, f(a)}$ such that the following three conditions hold.
\begin{enumerate}
\item Locally at $a$, $D \otimes_K \overline{K}$ is the locus of $\, \prod_i x_i = 0$.
\item Locally at $f(a)$, $E \otimes_K \overline{K}$ is the locus of $\, \prod_j y_j = 0$.
\item The morphism $f$ gives the $y_j$ as monomials in the $x_i$.
\end{enumerate}
Here we say that elements $z_1, \dots, z_n$ of a local ring $A$, containing its residue field, are \emph{uniformizing parameters} for $A$ if these elements minus their images in the residue field, form a system of regular parameters for $A$.
\end{named}

\begin{remark} \label{RemarkToroidal1} \rm
Note that, using the notation in the above definition, if $f$ is toroidal with respect to $D$ and $E$, then $f$ is log-smooth  with respect to $D$ and $E$ at each point of $X$.
The converse is also true, by the work of Kazuya Kato on logarithmic geometry (but we will not use this fact in the present paper).
\end{remark}

The following theorem is a small extension, proved in \cite{Abramovich-Denef-Karu}, of the Weak Toroidalization Theorem of Abramovich and Karu \cite{Abramovich-Karu}.

\begin{named}{Weak Toroidalization Theorem}
\label{TheoremWeakToroidalization}
Let $K$ be a field of characteristic zero. Let $f: X\to Y$ be a dominant morphism of varieties over $K$, and let $Z\subset X$ be a proper closed subset. 
Then there exist nonsingular quasi-projective varieties $X'$, $Y'$ over $K$,  and a commutative diagram
$$
\begin{array}{lclcl} 
D & \subset & X' &\stackrel{m_X}{\to} &X \\
& & \downarrow f' & & \downarrow f\\ 
E & \subset & Y' & \stackrel{m_Y}{\to} &Y
\end{array}
$$
with $m_X$, $m_Y$ projective birational morphisms over $K$, and $D, E$ reduced strict normal crossings divisors over $K$, such that

\begin{enumerate}
\item $f'$ is a quasi-projective morphism over $K$ and is toroidal with respect to $D$ and $E$,
\item $m_X^{-1}(Z)$ is a divisor on $X'$, and is contained in $D$,
\item the restriction of the morphism $m_X$ to $X'\setminus D$ is an open embedding.
\end{enumerate}
\end{named}
In the present paper, assertion (3) in the above theorem will not be used. The theorem is very much related to (but not implied by) Cutkosky's Theorem on Local Monomialization of Morphisms (Theorem 1.3 in \cite{Cutkosky-LocMon}). It is conjectured that we can take $m_X$ and $m_Y$ to be compositions of blow-ups of non-singular subvarieties. This conjecture is a weakening of the Conjecture on (Strong) Toroidalization \cite{Abramovich-Karu},\cite{Cutkosky-Tor3folds}.
Finally we mention that recently Illusie and Temkin obtained a result (Theorem 3.9 of \cite{Illusie-Temkin}) which is more general than the above Theorem \ref{TheoremWeakToroidalization}, and that Cutkosky \cite{Cutkosky-AnalyticMon} extended his Local Monomialization Theorem to complex and real analytic maps. 
% This might lead to an alternative approach to subanalytic sets. 

\section{The Tameness Theorem}  \label{SectionTamenessTheorem}

Let $R$ be a noetherian integral domain. In this section we assume that $R$ has characteristic zero. 
Let $f : X \rightarrow Y$ be a morphism of varieties over $R$.

\begin{named}{Tameness Theorem} \label{TheoremTameness}
Given rational functions $x_1, \dots, x_r$ on $X$, there exist rational functions $y_1, \dots, y_s$ on $Y$, and $\Delta \in R\setminus \{0\}$, such that for any $R[\Delta^{-1}]$-algebra A which is a henselian valuation ring we have the following.
Any $b \in Y(A)$ having the same residues with respect to $\, y_1, \dots, y_s$ as an image $f(a_0)$, with $a_0 \in X(A)$, is itself an image of an $a \in X(A)$ with the same residues as $a_0$ with respect to $\, x_1, \dots, x_r$.
\end{named}

\begin{remark} \label{RemarkAboutTameness1} \rm
In the statement of the Tameness Theorem we can choose the rational functions $y_1, \dots, y_s$ to be regular on $Y$, if $Y$ is affine. Indeed, this is a direct consequence of Lemma \ref{LemmaResidues1}.
\end{remark}

\subsection{Proof of the Tameness Theorem}
The Tameness Theorem can be proved easily by using Basarab's Elimination of Quantifiers \cite{Basarab-1}. Basarab's work is based on model theory using the same methods as Ax-Kochen and Er{\v{s}}ov. In the next subsections (\ref{ProofTamenessPart1}) up to (\ref{ProofTamenessPart5}) we present a purely algebraic geometric proof of the Tameness Theorem.

\subsubsection{Preliminary reductions.} \label{ProofTamenessPart1}
Let $K$ be the field of fractions of $R$. Our proof of the Tameness Theorem is by induction on the dimension of $X\otimes_RK$.

Covering $Y$ with a finite number of affine open subschemes, we see that in order to prove the Tameness Theorem we may suppose that $Y$ is affine. Moreover we may also suppose that $f$ is dominant. Indeed assume that $Y$ is affine and let $\overline{f(X)}$ be the Zariski closure of $f(X)$. Assume that the Tameness Theorem for the dominant morphism 
$X \rightarrow \overline {f(X)}$, induced by $f$, holds for given rational functions $x_1, \dots, x_r$ on $X$ and suitable chosen rational functions $y_1, \dots, y_s$ on $\overline{f(X)}$. By Lemma \ref{LemmaResidues1} we can actually choose the $y_1, \dots, y_s$ so that they are regular rational functions on $\overline{f(X)}$. Then obviously the Tameness Theorem for $f: X \rightarrow Y$ holds for the given $x_1, \dots, x_r$ and any finite list of regular rational functions on $Y$ which contains an extension to $Y$ for each of the regular rational functions $y_1, \dots, y_s$ on $\overline{f(X)}$, and which contains a sequence of functions whose zero locus equals $\overline{f(X)}$. It is clear that such a finite list exists because $Y$ is affine.

Covering $X$ with a finite number of affine open subschemes, and using Lemma \ref{LemmaResidues1}, we can further assume that $X$ is affine and that $x_1, \dots, x_r$ are regular rational functions on $X$.

Finally we can suppose that $R$ is a normal ring, because any finitely generated subring of $R$ becomes normal after inverting a suitable non-zero element (see e.g. section 32 of Matsumura \cite{Matsumura}).

\subsubsection{Applying the Weak Toroidalization Theorem.}
\label{ProofTamenessPart2}
Thus we assume that $R$ is normal, that $X$ and $Y$ are affine, that $f$ is dominant, and that  $x_1, \dots, x_r$ are regular rational functions on $X$. Moreover we suppose that no one of the given rational functions $x_i$ on $X$ is identically zero, because we can discard those that are identically zero. Let $Z$ be the union of the zero loci of the $x_i$ for $i = 1, \dots, r$.

Applying the Weak Toroidalization Theorem \ref{TheoremWeakToroidalization} and Remark \ref{RemarkToroidal1}, to the base change over $K$ of the above $X$, $Y$, $f$, and $Z$, we obtain a suitable $\Delta \in R \setminus \{0\}$, smooth quasi-projective varieties $X'$ and $Y'$ over $R[\Delta ^{-1}]$, and a commutative diagram of morphisms over $R[\Delta ^{-1}]$
$$
\begin{array}{lclccc} 
	D & \subset & X' &\stackrel{m_X}{\to} &X_0 & :=  X\otimes_R R[\Delta ^{-1}] \\
    & & \downarrow f' & & \downarrow f & \\ 
    E & \subset & Y'  & \stackrel{m_Y}{\to} &Y_0 & := Y\otimes_R R[\Delta ^{-1}] \\
\end{array}
$$
with $m_X$, $m_Y$ projective birational morphisms, $D \subset X'$, $E \subset Y'$ reduced strict normal crossings divisors over 
$R[\Delta ^{-1}]$ with $f^{-1}(E) \subset D$, such that
the following two conditions hold.
\begin{enumerate}
\item $f'$ is log-smooth with respect to $D$, $E$, at each point of $X'$.
\item $m_X^{-1}(Z\otimes_R R[\Delta ^{-1}])$ is a divisor on $X'$ and is contained in $D$.
\end{enumerate}

Choose a nonempty open subscheme $V \subset Y_0$ such that the rational map $m_Y^{-1}$ is regular on $V$. Moreover choose a nonempty open subscheme $U \subset X_0$ such that the rational map $m_X^{-1}$ is regular on $U$, and such that $U$ is disjoint from $m_X(D)$ and $f(U) \subset V$.

\subsubsection{Applying the Logarithmic Hensel's Lemma.}
\label{ProofTamenessPart3}
Choose rational functions $y_1, \dots, y_s$ on $Y$ such that the following two conditions hold.
\begin{enumerate}
\item
For each point $P'$ of $Y'$ and for each irreducible component $E_j$ of $E$, at least one of the elements 
$m_Y^*(y_1), \dots$, or $m_Y^*(y_s)$ belongs to $\cO_{Y',P'}$ and generates the ideal of $E_j$ in $\cO_{Y',P'}$.
\item
For any field $k$ which is an $R[\Delta ^{-1}]$-algebra, and any two distinct points $c'_1, c'_2 \in Y'(k)$, there is a 
$j \in \{ 1, \dots, s \}$ such that $y_j \in \cO_{Y',c'_1} \cap \cO_{Y',c'_2}$ and $y_j(c'_1) \neq y_j(c'_2)$.
\end{enumerate}
One easily verifies that condition (2) can be satisfied using that $Y'$ is quasi-projective over $R[\Delta ^{-1}]$.

We claim that the Logarithmic Hensel's Lemma for the morphism $f': X' \rightarrow Y'$ implies that the Tameness Theorem for the morphism $f: X \rightarrow Y$ is true for the given regular rational functions $x_1, \dots, x_r$ and the chosen rational functions $y_1, \dots, y_s$, under the additional assumptions that $a_0(\eta_A) \in U$ and  $b(\eta_A) \in V$. Here, as before, $\eta_A$ denotes the generic point of $\Spec(A)$. 

To verify this claim we argue as follows. Note that $b$ lifts to a $b' \in Y'(A)$ with $m_Y(b')=b$, and that $a_0$ lifts to a $a_0' \in X'(A)$ with $m_X(a_0')=a_0$, because $A$ is a valuation ring, $m_Y$ and $m_X$ are proper morphisms, and the rational maps $m_Y^{-1}$, and $m_X^{-1}$ are regular at respectively $b(\eta_A)$, and $a_0(\eta_A)$. Whence $a_0' \notin D(A)$,  since $U$ is disjoint from $m_X(D)$. Note also that $m_Y^{-1}$ is regular at $f(a_0)(\eta_A)$, since $f(U) \subset V$. Hence any rational function on $Y'$ which is regular at $f'(a_0')(\eta_A)$, respectively $b'(\eta_A)$, is also regular at  $f(a_0)(\eta_A)$, respectively $b(\eta_A)$, considered as rational function on $Y$, since $m_Y(f'(a_0')) = f(a_0)$.
Denoting, as before, the maximal ideal of $A$ by $\frm_A$, this implies that $b'(\frm_A) = f'(a_0')(\frm_A)$, by condition (2) with $k = A/\frm_A$, and the hypothesis that $b$ and $f(a_0)$ have the same residues with respect to $y_1, \dots, y_s$. Hence, using condition (1), we obtain that $b'$ and $f'(a_0')$ have the same residues with respect to the irreducible components of $E$.
\newline
We can now apply the Logarithmic Hensel's Lemma to the log-smooth morphism $f': X' \rightarrow Y'$, to get a point $a' \in X'(A)$ with $f'(a')=b'$, having the same residues as $a_0'$ with respect to the irreducible components of $D$. Hence $a'$ and $a_0'$ have the same residues with respect to $m_X^*(x_1), \dots, m_X^*(x_r)$, because of condition (2) in subsection \ref{ProofTamenessPart2} and Lemma \ref{LemmaResidues2}. Set $a := m_X(a')$. It is now obvious that $a$ satisfies the conclusion of the Tameness Theorem, because $x_1, \dots, x_r$ are regular.
\newline
Thus we have now verified our claim that the Tameness Theorem for the morphism $f: X \rightarrow Y$ is true for the given regular rational functions $x_1, \dots, x_r$ and the chosen rational functions $y_1, \dots, y_s$, under the additional assumptions that $a_0(\eta_A) \in U$ and  $b(\eta_A) \in V$.

Next we enlarge the list of rational functions $y_1, \dots, y_s$ on $Y$ by adjoining a sequence of regular rational functions on $Y_0$ whose zero locus on $Y_0$ equals $Y_0 \setminus V$. This is possible since $Y_0$ is affine. Then the condition $b(\eta_A) \in V$ is automatically satisfied if $a_0(\eta_A) \in U$, because then $f(a_0)(\eta_A) \in V$, since $f(U) \subset V$, and because $b$ has the same residues as $f(a_0)$ with respect to $y_1, \dots, y_s$.
\newline
Thus we have now proved the Tameness Theorem for the given regular rational functions $x_1, \dots, x_r$ and the chosen rational functions $y_1, \dots, y_s$, under the additional assumption that $a_0(\eta_A) \in U$. It remains to treat the case $a_0 \in (X_0 \setminus U)(A)$.

\subsubsection{Using the induction hypothesis.}
\label{ProofTamenessPart4}
Let $S$ be an irreducible component of $X_0 \setminus U$. Note that $\Dim(S \otimes_R K) < \Dim(Y \otimes_R K)$.
Let $s_i$ be the restriction to S of the regular rational function $x_i$, for $i=1, \dots, r$.
By the induction hypothesis the Tameness Theorem is true for the morphism $S \rightarrow Y$, induced by $f$, and the list of rational functions $s_1, \dots, s_r$. This yields a list of rational functions on $Y$ satisfying the Tameness Theorem for the restriction of $f$ to $S$. 

\subsubsection{Conclusion of the proof of the Tameness Theorem.}
\label{ProofTamenessPart5}
Replace the list $y_1, \dots, y_s$ by its union with the lists of rational functions on $Y$ obtained as above for each irreducible component $S$ of $X_0 \setminus U$. Then it is clear that the Tameness Theorem for $f: X \rightarrow Y$ holds for the given $x_1, \dots, x_r$ and the new list $y_1, \dots, y_s$. This finishes the proof of the Tameness Theorem. $\square$

\begin{remark} \label{RemarkAboutTameness2} \rm
To prove the Tameness Theorem in the special case that the list of rational functions $x_1, \dots, x_r$ on $X$ is empty, we can use Cutkosky's Local Monomialization Theorem \cite{Cutkosky-LocMon} instead of weak toroidalization. This special case is sufficient to prove Theorems \ref{TheoremTransferSurjectivity} and \ref{TheoremAxKochenForms}. 
Very recently Cutkosky \cite{Cutkosky-StrongerMonomialization} proved a stronger version of his Local Monomialization Theorem (with an additional requirement similar to (2) in \ref{TheoremWeakToroidalization}) that suffices to prove the Tameness Theorem in general.
\end{remark}

\section{Transfer of residues} \label{SectionTransferResidues}
In the next lemma we use the notation of the beginning of section \ref{SectionMultRes} and the terminology of subsection \ref{SubSecStructureMR}. Our proof of this lemma is an easy application of embedded resolution of singularities and does not depend on the Weak Toroidalization Theorem or the Tameness Theorem.

\begin{lemma} \textbf{\emph{Transfer of residues.}}
\label{LemmaTransferResidues}
Let $X$ be a variety over $\Z$, and $x_1, \dots, x_r$ regular rational functions on $X$. Then there exists a positive integer N such that for any two henselian valuation rings $A$ and $B$, having residue field characteristic $> N$ or zero, and admitting an isomorphism $\tau: \cM\cR(A) \rightarrow \cM\cR(B)$, we have the following. For any $a \in X(A)$ there exists $b \in X(B)$ such that $a \!\mod \frm_A \; = \; b \!\mod \frm_B$, identifying $A / \frm_A$ with $B / \frm_B$ via $\tau$, and such that $\tau(\mres (x_i(a))) = \mres (x_i(b))$, for $i=1, \cdots, r$. 
\end{lemma}

\emph{Proof}. Let $Z$ be the union of the zero loci of the regular rational functions $x_1, \dots, x_r$ on $X$. By using embedded resolution of $Z \otimes \Q \subset X \otimes \Q$ and induction on $\dim X$, and by inverting finitely many primes, we can assume that $X$ is smooth over $\Z$ and that $Z$ is a strict normal crossings divisor over $\Z$. 
Note that this reduction requires $A$ and $B$ to be valuation rings in order to apply the valuative criterion of properness to the resolution morphism.
Moreover, covering $X$ with finitely many suitable open subschemes, we can further assume that there exist uniformizing parameters $z_1,\dots, z_n$ over $\Z$ on $X$ such that $Z$ is the locus of $\prod_{j=1, \dots, n} z_j$ and such that the locus of each $z_i$ is irreducible or empty. Then each $x_i$ is a monomial in the $z_j$'s times a unit in $\Gamma(X, \cO_X$), because $X$ is normal. Hence it suffices to prove the lemma for $x_1, \dots, x_r$ replaced by the uniformizing parameters $z_1,\dots, z_n$. Thus we can assume that $x_1, \dots, x_r$ are uniformizing parameters over $\Z$ on $X$. But then the lemma is a direct consequence of Hensel's Lemma for the \'{e}tale morphism $\pi: X \rightarrow \A^r$ induced by $x_1,\dots, x_r$. Indeed, for $i=1, \dots, r$, choose $\beta_i \in B$ such that $\mres(\beta_i) = \tau(\mres (x_i(a)))$. Then, by Hensel's Lemma, there exists $b \in X(B)$ such that $\pi(b) = (\beta_1, \dots, \beta_r)$ and  $a \!\mod \frm_A \; = \; b \!\mod \frm_B$, identifying $A / \frm_A$ with $B / \frm_B$ via $\tau$. This rational point  $b$ satisfies the requirements of the lemma. $\square$

\section{Transfer of surjectivity}  \label{SectionTransferSurjectivity}
In this section we prove Theorem \ref{TheoremTransferSurjectivity} on transfer of surjectivity. We use the terminology of subsection \ref{SubSecStructureMR}. For any prime number $p$ there is an obvious unique isomorphism $\tau_p: \cM\cR(\Z_p) \rightarrow \cM\cR(\F_p[[t]])$ such that, for any integer $m$ and any unit $u$ in $\Z_p$, we have $\tau_p(\mres(up^m)) = \mres((u \!\mod p)t^m)$.

\begin{remark} \label{RemarkIdentifyingResidues} \rm
By using the isomorphism $\tau_p$ to identify multiplicative residues of elements of $\Z_p$ with multiplicative residues of elements of $\F_p[[t]]$, one can give an obvious meaning to Definition \ref{DefinitionMultRes1}, with $X$ a variety over $\Z$, when $a \in X(\Z_p)$ and $a' \in X(\F_p[[t]])$. Thus it is well defined  when we say that $a$ and $a'$ \emph{have the same residues with respect to} given rational functions $x_1, \dots, x_r$ on $X$. We will use this terminology throughout the following subsection.
\end{remark}

\subsection{Proof of Theorem \ref{TheoremTransferSurjectivity} on transfer of surjectivity.}

Covering $Y$ with a finite number of affine open subschemes, we can assume that $Y$ is affine. By the Tameness Theorem \ref{TheoremTameness}, with $R = \Z$, and by Remark \ref{RemarkAboutTameness1}, there exist regular rational functions $y_1, \dots, y_s$ on $Y$, and $\Delta \in \Z \setminus \{0\}$, such that for any prime $p$, which does not divide $\Delta$, the following is true for both $A = \Z_p$ and $A = \F_p[[t]]$. Any $b \in Y(A)$ having the same residues with respect to $\, y_1, \dots, y_s$ as an image $f(a_0)$, with $a_0 \in X(A)$, is itself an image of an $a \in X(A)$.

Let $p$ be any prime which does not divide $\Delta$, and which is large enough as required by the statement of Lemma \ref{LemmaTransferResidues}, on transfer of residues, for both the list of regular rational functions $y_1, \dots, y_s$ on $Y$ and the list of regular rational functions $y_1 \circ f, \dots, y_s \circ f$ on $X$.

Assume now that $f: X(\F_p[[t]]) \rightarrow Y(\F_p[[t]])$ is surjective. We are going to prove that $f: X(\Z_p) \rightarrow Y(\Z_p)$ is surjective. The implication in the other direction is proved in exactly the same way.

Take any $b \in Y(\Z_p)$. By transfer of residues there exists a point $b' \in Y(\F_p[[t]])$ having the same residues as $b$ with respect to $y_1, \dots, y_s$, in the sense of Remark \ref{RemarkIdentifyingResidues}. By the surjectivity assumption, there exists a point $a' \in X(\F_p[[t]])$ with $f(a') = b'$. By transfer of residues, we find a point $a_0 \in X(\Z_p)$ having the same residues as $a'$ with respect to  $y_1 \circ f, \dots, y_s \circ f$. Notice that, by construction, $f(a_0)$, $f(a')$, and $b$ have the same residues with respect to $y_1, \dots, y_s$. Hence, by the above mentioned instance of the Tameness Theorem, the point $b$ is the image of a rational point $a \in X(Z_p)$. This proves the surjectivity of $f: X(\Z_p) \rightarrow Y(\Z_p)$. $\square$

\section{Elimination of quantifiers} \label{SectionEliminationQuantifiers}

\subsection{The languages $\LRings$ and $\LMR$} \label{SubsectionLanguages}
A \emph{formula} in the \emph{language $\LRings$ of rings}  is built from the \emph{logical connectives} (and, or, not, implies, iff), variables, existential and universal quantifiers, equality, addition, multiplication, 0, and 1.
An \emph{assertion} (sentence) in $\LRings$ is a formula in $\LRings$ without free variables. Using formulas in the language $\LRings$ we can talk about any ring $A$, interpreting the variables and the quantifiers as running over $A$.

A \emph{formula in the language $\LMR$ of multiplicative residues} is built from the logical connectives, variables, existential and universal quantifiers,  equality, multiplication, the binary composition law $+_{\mathrm{mod}}$, 0, and 1. Using formulas in the language $\LMR$ we can talk about $\cM\cR(A)$ for any local integral domain $A$, interpreting the variables and the quantifiers as running over $\cM\cR(A)$, and interpreting $+_{\mathrm{mod}}$ by the binary composition law $+_{\mathrm{mod}}$ introduced in subsection \ref{SubSecStructureMR}.

A \emph{formula in the language} $\LRingsMR$ \emph{of rings with multiplicative residues} is built from the logical connectives, variables called of the \emph{first sort}, variables called of the \emph{second sort}, existential and universal quantifiers with respect to variables of the first or second sort, formulas in $\LRings$, with variables of the first sort, and expressions obtained from formulas in $\LMR$, with variables of the second sort, by replacing some (or none) of the free variables by the operator $\mres$ applied to polynomials over $\Z$ in variables of the first sort. 
Using formulas in the language $\LRingsMR$ we can talk about any local integral domain $A$, interpreting the variables of the first sort (and quantifiers with respect to these) as running over $A$, and interpreting the variables of the second sort (and quantifiers with respect to these) as running over $\cM\cR(A)$. 
Note that each formula in $\LRings$ and each formula in $\LMR$ is also a formula in $\LRingsMR$.

Below, we use the symbol $\wedge$ to denote the logical connective \textquotedblleft and\textquotedblright.

We recall that for any local integral domain $A$ we denote by $s_A$ the  natural map $\cM\cR(A) \rightarrow A / \frm_A$  introduced in subsection \ref{SubSecStructureMR}.

\begin{lemma} \textbf{\emph{Residue field interpretability in $\LMR$.}}
\label{LemmaModularAddition}
\newline
Let $\varphi(x_1, \dots, x_n)$ be a formula in the language $\LRings$, with free variables $x_1, \dots, x_n$. Then there exists a formula $\theta(\xi_1, \dots, \xi_n)$ in the language $\LMR$, with free variables $\xi_1, \dots, \xi_n$, such that for any local integral domain $A$ and any $a_1, \dots, a_n \in \cM\cR(A)$ we have the following.  The formula $\varphi(s_A(a_1), \dots, s_A(a_n))$ holds in $A/\frm_A$ if and only if the formula $\theta(a_1, \dots, a_n)$ holds in $\cM\cR(A)$.  
\end{lemma}
\emph{Proof}.
This is straightforward and left to the reader, using the following obvious observations for any $a_1, a_2, a_3 \in \cM\cR(A)$.
\begin{enumerate}
\item
$s_A(a_1) = 0$ if and only if $\, a_1 +_{\mathrm{mod}} 0 = 0$.
% not $\exists \, c \in \cM\cR(A): a_1\,c = 1$.
\item 
$s_A(a_1) = s_A(a_2)$ if and only if ($a_1 = a_2$ or $0 = s_A(a_1) = s_A(a_2)$).
\item
$s_A(a_1)+s_A(a_2)=s_A(a_3)$ if and only if $s_A(a_1 +_{\mathrm{mod}} a_2) = s_A(a_3)$.
\item
$s_A(a_1)s_A(a_2)=s_A(a_3)$ if and only if $s_A(a_1 a_2)=s_A(a_3)$.
\item
$\exists \, x \in A/\frm_A \!: S(x)$ if and only if $\, \exists \, y \in \cM\cR(A) \!: S(s_A(y))$, where $S$ is any relation in one variable on $A/\frm_A$. $\square$
\end{enumerate}

\begin{lemma} \textbf{\emph{Realizability of multiplicative residues.}} \label{Lemma-RealizabilityOfResidues}
Let $X$ be a variety over $\Z$ and $x_1, \dots x_r$ regular rational functions on $X$. Then there exists a formula $\psi(\lambda_1, \dots, \lambda_r)$ in the language $\LMR$, with free variables $\lambda_1, \dots, \lambda_r$, and a positive integer $N$, such that for any henselian valuation ring $A$, having residue field characteristic $> N$ or zero, and for any $b_1, \dots, b_r \in \cM\cR(A)$ we have  the following. There exists an $a \in X(A)$ with $\mres(x_i(a) )=b_i$, for each $i=1, \dots, r$, if and only if $\psi(b_1, \dots, b_r)$ holds in $\cM\cR(A)$ . 
\end{lemma}

\emph{Proof}. Our proof is an easy application of embedded resolution of singularities and does not depend on the Weak Toroidalization Theorem or the Tameness Theorem. Reasoning as in the proof of Lemma \ref{LemmaTransferResidues} (based on embedded resolution of singularities) we can assume that $X$ is an affine smooth variety over $\Z$, and that there exist uniformizing parameters $z_1, \dots, z_n$ over $\Z$ on $X$, such that each $x_i$, for $i=1, \dots, r$, can be written as a monomial in $z_1, \dots, z_n$ times a unit $u_i$ in $\Gamma(X, \cO_X)$, i.e.
\begin{equation} \label{equation-RealizabilityResidues-1}
\nonumber
x_i=u_i\prod_{j=1}^n z_j^{e_{i,j}},
\end{equation}
for i=1, \dots, r, where the $e_{i,j}$ are nonnegative integers.
Hence the relation
\begin{equation} \label{equation-RealizabilityResidues-2}
\exists \, a \in X(A) \; \forall i \in \{1, \dots, r\}:
\mres(x_i(a) )=b_i
\end{equation}
is equivalent to
\begin{eqnarray}
\nonumber
& \exists \, a \in X(A) \; 
\exists \, c_1, \dots, c_n, c'_1, \dots, c'_r \in \cM\cR(A)\!: 
\bigwedge_{j=1}^n \;
\bigwedge_{i=1}^r \\
\nonumber
& c_j =  \mres(z_j(a)), \;\; c'_i = \mres(u_i(a)), \;\; 
c'_i \prod_{j=1}^n c_j^{e_{i,j}}=b_i .
\end{eqnarray}
Because, for any $a \in X(A)$, $u_i(a)$ is a unit in $A$, the condition $c'_i = \mres(u_i(a))$ is equivalent to
$s_A(c'_i) = u_i(a) \!\! \mod \frm_A$.
Thus the relation (\ref{equation-RealizabilityResidues-2}) is equivalent to
\begin{eqnarray}
\nonumber
& \exists \, a \in X(A) \; 
\exists \, c_1, \dots, c_n, c'_1, \dots, c'_r \in \cM\cR(A)\!: 
\bigwedge_{j=1}^n \;
\bigwedge_{i=1}^r \\
\nonumber
& c_j =  \mres(z_j(a)), \; \;
s_A(c'_i) = u_i(a) \!\!\! \mod \frm_A, \; \;
c'_i \prod_{j=1}^n c_j^{e_{i,j}}=b_i .
\end{eqnarray}
Hence, applying Hensel's Lemma for the \'{e}tale morphism $\pi: X \rightarrow \A^n$ induced by $z_1,\dots, z_n$, we see that the relation (\ref{equation-RealizabilityResidues-2}) is equivalent to
\begin{eqnarray}
\nonumber
& \exists \, \bar{a} \in X(A/\frm_A) \; 
\exists \, c_1, \dots, c_n, c'_1, \dots, c'_r \in \cM\cR(A)\!: 
\bigwedge_{j=1}^n \;
\bigwedge_{i=1}^r \\
\nonumber
& s_A(c_j) =  z_j(\bar{a}), \; 
s_A(c'_i) = u_i(\bar{a}), \; 
c'_i \prod_{j=1}^n c_j^{e_{i,j}}=b_i .
\end{eqnarray}
Indeed, by Hensel's Lemma, any $\bar{a} \in X(A/\frm_A)$ can be lifted to an $a \in X(A)$ with $z_1(a) = \gamma_1, \dots, z_n(a) = \gamma_n$, for any $\gamma_1, \dots, \gamma_n \in A$ satisfying $z_j(\bar{a}) = \gamma_j \! \mod \frm_A$, for $j=1, \dots, n$. We conclude that the relation (\ref{equation-RealizabilityResidues-2}) can be expressed by a formula in the language $\LMR$, because of Lemma \ref{LemmaModularAddition}. 
This terminates the proof of the lemma. $\square$

\subsection*{}
The next result is a reformulation of Basarab's Quantifier Elimination Theorem \cite{Basarab-1}, for henselian valuation rings with residue field characteristic large enough or zero. It eliminates the quantifiers running over the valuation ring at the expense of introducing new quantifiers running over the multiplicative residues (which  are often easier to analyze). Our proof of Basarab's result is an easy application of the Tameness Theorem \ref{TheoremTameness} and the previous lemma.

\begin{named}{Basarab's Quantifier Elimination Theorem}
\label{TheoremBasarab}
Let $\varphi(x, \gamma)$ be a formula in the language $\LRingsMR$, with free variables $x =(x_1, \dots, x_n)$ of the first sort and free variables $\gamma = (\gamma_1, \dots, \gamma_r)$ of the second sort. Then there exist a positive integer $N$, polynomials $P_1(x), \dots, P_m(x)$ over $\Z$, and a formula $\theta(\xi_1, \dots, \xi_m, \gamma)$ in the language $\LMR$, with free variables $\xi_1, \dots, \xi_m, \gamma$, such that the equivalence
$$
\varphi(x, \gamma) \leftrightarrow \theta(\mres(P_1(x)), \dots, \mres(P_m(x)), \gamma)
$$
holds for any henselian valuation ring having residue field characteristic larger than $N$ or zero. 
In particular, this holds for any formula $\varphi$ in $\LRings$; then there is no $\gamma$ involved.

\end{named}

\emph{Proof}. We can assume that the formula $\varphi(x, \gamma)$ does not contain universal quantifiers, because a universal quantifier can be expressed by the negation of an existential quantifier.  By induction on the sum of the number of quantifiers and the number of logical connectives in $\varphi(x, \gamma)$, we can suppose that $\varphi(x, \gamma)$ starts with an existential quantifier, with respect to a variable of the first sort, and is equivalent to
\begin{equation} \label{equation-BasarabQE-1}
\exists \, y_1: \psi(\mres(P_1(x, y_1)), \dots, \mres(P_k(x, y_1)), \gamma),
\end{equation}
with $P_1(x, y_1), \dots, P_k(x, y_1) \in \Z[x,y_1]$ and $\psi$ a formula in $\LMR$. Here and below, with \textquotedblleft equivalent\textquotedblright  we mean equivalent for all henselian valuation rings with residue field characteristic large enough or zero. Obviously, (\ref{equation-BasarabQE-1}) is equivalent to
\begin{equation} \label{equation-BasarabQE-2}
\exists \, \lambda_1, \dots, \lambda_k: \left(  \psi( \lambda_1, \dots, \lambda_k, \gamma) \wedge \exists \, y_1: \bigwedge_{j=1}^k \mres (P_j(x, y_1)) = \lambda_j \right).
\end{equation}

Applying the Tameness Theorem \ref{TheoremTameness} and Remark \ref{RemarkAboutTameness1} to the morphism $\A^{n+1} \rightarrow \A^n \! : (x,y_1) \mapsto x$ and the functions $P_1(x, y_1), \dots, P_k(x, y_1)$ on  $\A^{n+1}$, yields functions $Q_1(x), \dots, Q_s(x) \in \Z[x]$ on $\A^n$ such that 
$$\exists \, y_1 \!: \bigwedge_{j=1}^k \mres (P_j(x, y_1)) = \lambda_j$$
is equivalent to
\begin{equation}
\nonumber
\exists \, x',y'_1 \!: \bigwedge_{j=1}^k \mres (P_j(x', y'_1)) = \lambda_j \wedge
\bigwedge_{l=1}^s \mres (Q_l(x')) = \mres (Q_l(x)).
\end{equation}
Here $x'=(x'_1, \dots, x'_n)$. This last equivalence implies the theorem, because $\varphi(x, \gamma)$ is equivalent to (\ref{equation-BasarabQE-2}) and because, by Lemma \ref{Lemma-RealizabilityOfResidues}, the last formula above is equivalent to a formula of $\LMR$ with some of its free variables replaced by $\mres(Q_1(x)), \dots, \mres(Q_s(x))$. 
This concludes the proof of the theorem. $\square$.

\subsection{Proof of the Ax-Kochen-Er{\v{s}}ov Transfer Principle \ref{TheoremAKE-TransferPrinciple}}
This is a direct consequence of Basarab's Theorem \ref{TheoremBasarab}, because there is an isomorphism $\tau_p: \cM\cR(\Z_p) \rightarrow \cM\cR(\F_p[[t]])$ as explained in the beginning of section \ref{SectionTransferSurjectivity}.  $\square$

\bibliographystyle{plain}
% \bibliography{JD}

\end{document}